\documentclass[12pt]{article}%
\usepackage{amsmath}
\usepackage{amsfonts}
\usepackage{amssymb}
\usepackage{graphicx}%
\newtheorem{theorem}{Theorem}

\newtheorem{condition}{Condition}

\newtheorem{lemma}{Lemma}

\newtheorem{proposition}{Proposition}
\newtheorem{remark}{Remark}

\begin{document}

\title{On maximum of Gaussian random field having unique maximum point of its
variance\thanks{Partially supported by Russian Science Foundation, grant
14-49-00079. The authors thank Enkelejd Hashorva for fruitful discussions.}}
\author{Sergey G. Kobelkov\thanks{Lomonosov Moscow State University, Moscow, Russia,
\emph{{sergeyko81@gmail.com}}}
\and Vladimir I. Piterbarg\thanks{Lomonosov Moscow state university, Moscow,
Russia; Scientific Research Institute of System Development of the Russian
Academy of sciences; Federal National Research University ``Moscow Power
Engineering Institute'' \emph{{piter@mech.math.msu.su}}}}
\date{}
\maketitle

\textbf{Abstract:} Gaussian random fields on Euclidean spaces whose variances
reach their maximum values at unique points are considered. Exact asymptotic
behaviors of probabilities of large absolute maximum of theirs trajectories
have been evaluated using Double Sum Method under the widest possible conditions.

\


\section{ Introduction.}

This contribution is a generalization of results of \cite{KHP}. As discussed
in \cite{PP}, and then in \cite{book}, \cite{lectures}, \cite{KHP},
non-stationary Gaussian processes are more subtle to deal with since both the
local properties of the variance function at its point of global maximum and
those of the covariance function have to be carefully formulated. One can say
the same about Gaussian fields, see \cite{book}. Our aim is to show maximum
capabilities of the Pickands' Double Sum Method applying to Gaussian fields
with unique point of global maximum; in case of processes this has been done
in \cite{KHP}. The Pickands' method was developed originally for asymptotic
behavior of the maximum tail distribution for Gaussian stationary processes in
\cite{pickands}, with corrections in \cite{pit1}. This method has been
generalized to Gaussian random fields, \cite{book}, where stationary fields
with power like behavior of the covariance function at zero are considered as
well as fields with similar behavior of the covariance function at the unique
maximum point of variance. However, while the power behavior of the covariance
function, with possible slight generalization to regular variation of it,
\cite{pit1}, is quite essential for the Pickand's method, the required in
\cite{PP}, \cite{book} power behavior of the variance, as it has been shown in
\cite{KHP}, looks somewhat artificial. Thus the principal task of this
contribution is to investigate the tail asymptotic behavior of supremum of
non-stationary Gaussian fields by imposing a weak and natural assumption on
their variance functions, see Conditions \ref{CondA5} and \ref{CondA6} below.
In connection with, notice that in the recent article \cite{HDL} it is proved
that in the non-stationary case the behavior of variance does not need to be
exactly power but may be just regularly varying. Here we do not assume even this.

Let $S\subset\mathbb{R}^{d}$ be the closure of a bounded open set containing
zero, and let $X(\mathbf{t}),$ $\mathbf{t}\in S,$ be a zero mean a.s.
continuous Gaussian random field with covariance function $R(\mathbf{s,t)=E}%
X(\mathbf{s})X(\mathbf{t})$; denote by $\sigma^{2}(\mathbf{t})=R(\mathbf{t,t}%
)$ its variance function, which is continuous since $X$ is a.s. continuous. We
study the asymptotic behavior of the probability%

\begin{equation}
P(S;u)=\mathbf{P}(\max_{\mathbf{t\in}S}X(\mathbf{t)>}u) \label{P(S,u)}%
\end{equation}
as $u\rightarrow\infty.$ We need a slightly stronger condition than a.s.
continuity of sample paths. Denote $\mathbb{B}_{\varepsilon}:=\{\mathbf{t}%
:|\mathbf{t}|\leq\varepsilon\}.$

\begin{condition}
\label{CondA1} $X(\mathbf{t})$ is a.s. continuous. Moreover, there exists
$\varepsilon>0$ such that Dudley's integral, \cite{Dudley}, \cite{lectures},
for the standardized field $\bar{X}(\mathbf{t})=X(\mathbf{t})/\sigma
(\mathbf{t})$, $\mathbf{t}\in\mathbb{B}_{\varepsilon}$ is finite.
\end{condition}

Notice that for homogeneous Gaussian fields this condition is also necessary
for existing of a.s. continuous version of the field (X. Fernique
\cite{fernique}). We need this condition in order to use V. A. Dmitrovsky's
inequality for estimating the exit probability from above. For reader's
convenience we give a corollary from the inequality adapted to our purposes,
see Corollary 8.2.1, \cite{lectures}.

\begin{proposition}
\label{Dmitrovsky}(V. Dmitrovsky, \cite{Dmitr1}, \cite{Dmitr2},
\cite{lectures}) Let Condition \ref{CondA1} be hold. Then there exists
$\gamma(u)$, such that $\gamma(u)\rightarrow0$ as $u\rightarrow\infty$, and
for any $S_{1}\subset\mathbb{B}_{\varepsilon}$,%
\[
P(S_{1};u)\leq\exp\left(  -\frac{u^{2}}{2\sigma^{2}(S_{1})}+u\gamma(u)\right)
,
\]
where $\sigma^{2}(S_{1})=\sup_{\mathbf{t}\in S_{1}}\mathbf{E}X^{2}%
(\mathbf{t})$
\end{proposition}

\begin{condition}
\label{CondA2} $\sigma(\mathbf{t})$ reaches its absolute maximum on $S$ at
only $\mathbf{0.}$
\end{condition}

Without loss of generality assume that $\sigma(\mathbf{0})=1.$ Notice also
that by time parameter shift the maximum point can be made arbitrary, with
corresponding conditions on the parameter set. From Condition \ref{CondA2} it
follows in particular that the normalized field $\bar{X}(\mathbf{t})$,
$\mathbf{t}\in\mathbb{B}_{\varepsilon},$ see Condition \ref{CondA1}, exists.
Furthermore, it follows from it that $r(\mathbf{s,t)}\leq1$ and the equality
holds only for $\mathbf{s=t=0.}$

Notice that we do not consider the case of boundary maximum point of variance.
It can be considered with described here tools, with involving the structure
of the boundary near the point. Such the consideration could not require any
new ideas but makes the text longer and even more difficult to read. We would
have to introduce a series new Pickands like constants. Only in the Talagrand
case (see below) the asymptotic behavior reminds the same in this boundary
maximum point case.

\begin{condition}
\label{CondA3} (Local stationarity at $\mathbf{0}$). There exists a covariance
function $r(\mathbf{t)}$ of a homogeneous random field with $r(\mathbf{t)}<1$
for all $\mathbf{t\neq0}$ such that
\[
\lim_{\mathbf{s,t\rightarrow0,s\neq t}}\frac{1-R(\mathbf{s,t)}}%
{1-r(\mathbf{t-s})}=1.
\]

\end{condition}

\begin{remark}
\label{cor_cov} In contrast to \cite{PP}, \cite{book}, \cite{lectures},
\cite{KHP}, and other works, we assume here local stationarity in terms of
covariance function, not correlation (normed covariance) function. This is
because we would like to impose minimal number of assumptions on the variance
function. In particular, we do not assume H\"older condition in any
neighborhood of zero, like we did this in \cite{KHP}.
\end{remark}

For vectors $\mathbf{a}=(a_{1},...,a_{d}),$ $\mathbf{b}=(b_{1},...,b_{d})$
define $\mathbf{ab}=(a_{1}b_{1},...,a_{d}b_{d}).$ For a set $T$ we write
$\mathbf{a}T=\{\mathbf{at,t\in}T\}.$

\begin{condition}
\label{CondA4} There exists a basis in $\mathbb{R}^{d}$, a vector function
$\mathbf{q}(u)=(q_{1}(u),...,q_{d}(u))$, $u>0$, $q_{i}(u)>0,\ i=1,...,d,$ and
a positive for all $\mathbf{t\neq0}$ function $h(\mathbf{t})$ such that for
any $\mathbf{t}$ written in these coordinates,
\begin{equation}
\lim_{u\rightarrow\infty}u^{2}(1-r(\mathbf{q}(u)\mathbf{t))=}h(\mathbf{t),}
\label{A3}%
\end{equation}
uniformly in $\mathbf{t}$ from any closed set.
\end{condition}

In slightly other words, Condition \ref{CondA4} means that for some orthogonal
matrix $U$, (\ref{A3}) is fulfilled for $\mathbf{t}^{\prime}=U\mathbf{t}$
instead of $\mathbf{t}$, that is%
\[
\lim_{u\rightarrow\infty}u^{2}(1-r(\mathbf{q}(u)U\mathbf{t))=}h(U\mathbf{t).}%
\]
Remark that by definition of uniform convergence to a positive for all
non-zero $\mathbf{t}$ function with $h(\mathbf{0})=0$, from Condition
\ref{CondA4} it follows that $h(\mathbf{t})$ is continuous, and for any
continuous function $c_{\mathbf{t}}$ with $\lim_{\mathbf{t\rightarrow0}%
}c_{\mathbf{t}}=1,$%
\begin{equation}
\lim_{\mathbf{t\rightarrow0}}\frac{1-r(c_{\mathbf{t}}\mathbf{t)}%
}{1-r(\mathbf{t})}=1\mathbf{.} \label{A4uniform}%
\end{equation}
Consider a simple example. Let $d=2$ and
\[
1-r(\mathbf{t})=(|t_{1}+t_{2}|^{\alpha_{1}}+|t_{1}-t_{2}|^{\alpha_{2}%
})(1+o(1))
\]
as $\mathbf{t\rightarrow0}$, where $2\geq\alpha_{1}>\alpha_{2}>0$. For such
covariance function one cannot find $\mathbf{q}(u)$ satisfying (\ref{A3}),
whereas rotating the basis turning at an angle of $\pi/4,$ we have,
\[
1-r(\mathbf{t})=(2^{\alpha_{1}/2}|t_{1}|^{\alpha_{1}}+2^{\alpha_{2}/2}%
|t_{2}|^{\alpha_{2}})(1+o(1))
\]
as $\mathbf{t\rightarrow0,}$ and $\mathbf{q}(u)=(u^{-2/\alpha_{1}%
},u^{-2/\alpha_{2}})$ with $h(\mathbf{t})=2^{\alpha_{1}/2}|t_{1}|^{\alpha_{1}%
}+2^{\alpha_{2}/2}|t_{2}|^{\alpha_{2}}.$

From now on we assume that the basis in $\mathbb{R}^{d}$ satisfies Condition
\ref{CondA4}. From Conditions \ref{CondA3} and \ref{CondA4} it follows a
regular variation property of $r(\mathbf{t})$. Indeed, remark first that since
$r(\mathbf{t})$ is continuous at zero, we have $\mathbf{q}(u)\rightarrow
\mathbf{0}$ as $u\rightarrow\infty.$ Denote $(\mathbf{e}_{i},i=1,...,d),$
coordinate vectors in $\mathbb{R}^{d}.$ From Condition \ref{CondA4}, since
$h(\mathbf{t)}$ is continuous, it also follows that
\[
\lim_{s\rightarrow0}\frac{1-r(st\mathbf{e}_{i})}{1-r(s\mathbf{e}_{i})}%
=\frac{h(t\mathbf{e}_{i})}{h(\mathbf{e}_{i})}\ \text{for all }t>0\text{ and
}i=1,...,d,
\]
and since the functions on the right are positive, continuous and cannot be
equal identically to one ($h(\mathbf{0})=0$)$,$ $1-r(t\mathbf{e}_{i})$
regularly varies at zero with positive degree $\alpha_{i}>0,$ $i=1,...,d,$ and
$h(t\mathbf{e}_{i})=C_{i}|t|^{\alpha_{i}},$ $i=1,...,d$ for some $C_{i}$.
Notice that from properties of positive defined functions it follows that
$\alpha_{i}\leq2$ for all $i.$ That is, denoting by $\ell_{i}(t),$
$i=1,...,d,$ the corresponding slowly varying functions, we write for all $i,$
that $1-r(t\mathbf{e}_{i})=\ell_{i}(t)|t|^{\alpha_{i}}.$ Hence, using
(\ref{A3}) and the definition of slowly variation,%
\[
u^{2}(1-r(q_{i}(u)t\mathbf{e}_{i}))=u^{2}q_{i}(u)^{\alpha_{i}}|t|^{\alpha_{i}%
}\ell_{i}(tq_{i}(u))=(1+o(1))u^{2}q_{i}(u)^{\alpha_{i}}|t|^{\alpha_{i}}%
\ell_{i}(q_{i}(u))
\]
as $u\rightarrow\infty,$ and
\[
u^{2}q_{i}(u)^{\alpha_{i}}\ell_{i}(q_{i}(u))\rightarrow1
\]
as $u\rightarrow\infty.$ From here it follows that for some $\ell_{1i}(u)$
\begin{equation}
q_{i}(u)=\ell_{1i}(u)u^{-2/\alpha_{i}}\text{ and }\ell_{1i}(u)^{\alpha_{i}%
}\ell_{i} (u^{-2/\alpha_{i}}\ell_{1i}(u))\rightarrow1\text{ as }%
u\rightarrow\infty. \label{qi}%
\end{equation}
It is shown in \cite{KHP} that
\begin{equation}
\ell_{1i}(u)=(1+o(1))\ell_{i}^{\#}(u^{-2})^{1/\alpha_{i}} \label{q(u)}%
\end{equation}
as $u\rightarrow\infty,$ where $\ell_{i}^{\#}$ is the de Bruijn conjugate of
$\ell_{i}$, see details below and also in \cite{BGT}, and $\ell_{1i}(u)$
slowly varies as well. Moreover, in \cite{KHP} it is shown that (\ref{A3})
holds for any $\ell_{i}^{\prime}(u),$ $i=1,...,d,$ such that $\lim
_{u\rightarrow\infty}\ell_{i}(u)/\ell_{i}^{\prime}(u)=1,$ $i=1,...,d.$
Consequently, without loss of generality we assume in the following that all
$\ell_{i}(u),$ $i=1,...,d$ are monotone. Furthermore, using this argument and
the fact that the ratio of two slowly varying positive functions, say
$\mathcal{L}_{1}(u)$ and $\mathcal{L}_{2}(u),$ slowly varies as well, having
in mind again its monotone equivalent, we may write that%
\begin{equation}
\lim_{u\rightarrow\infty}\mathcal{L}_{1}(u)/\mathcal{L}_{2}(u)=\kappa
\in\lbrack0,\infty]. \label{eqviv_L}%
\end{equation}
Let us say that $\mathcal{L}_{1}(u)\succ\mathcal{L}_{2}(u)$ if $\kappa
=\infty,$ $\mathcal{L}_{1}(u)\prec\mathcal{L}_{2}(u),$ if $\kappa=0,$ and
$\mathcal{L}_{1}(u)\asymp\mathcal{L}_{2}(u)$ if $\kappa\in(0,\infty).$ We use
these definitions in the proof of the following proposition.

\begin{proposition}
\label{e_reg_var}Let Conditions \ref{CondA3} and \ref{CondA4} be fulfilled for
a covariance function $r(\mathbf{t}).$ Then for any vector $\mathbf{f}$ the
function $r(t\mathbf{f)}$ regularly varies at zero with degree $\alpha
(\mathbf{f})\in(0,2].$ Moreover, if $\alpha(\mathbf{f})=2,$ then
$\lim_{t\rightarrow0}t^{-2}(1-r(t\mathbf{f})\dot{)}\in(0,\infty].$
\end{proposition}

\textbf{Proof: } Denote $g(s):=1-r(s\mathbf{f}),$ $\mathbf{f=(}f_{i}%
,i=1,...,d),$
\begin{equation}
\alpha=\alpha(\mathbf{f)}:=\min\{\alpha_{i}:f_{i}\neq0\},\text{ and
}\mathcal{J=J(}\mathbf{f}):=\{i:\alpha_{i}=\alpha,f_{i}\neq0\}.\label{alpha}%
\end{equation}
Introduce the \textquotedblleft main\textquotedblright\ index set
$\mathcal{J}_{0}\subset\mathcal{J}$ related to $\mathbf{f,}$ with property
\[
\forall i,j\in\mathcal{J}_{0},\ell_{1i}(u)\asymp\ell_{1j}(u),\text{
and\ }\forall i,j:i\in\mathcal{J}_{0},\ j\in\mathcal{J}\setminus
\mathcal{J}_{0},\ \ell_{1i}(u)\succ\ell_{1j}(u).
\]
Numerate indexes from $\mathcal{J}_{0}=\{i_{1},...,i_{k_{0}}\}$ and denote
\[
\kappa_{i}=\lim_{u\rightarrow\infty}\ell_{i}(u)/\ell_{1i_{1}}(u)\in
(0,\infty),\ i=i_{1},...,i_{k_{0}},
\]
where, as above, we mean monotone equivalent ratios, see (\ref{eqviv_L}). Put
$s=u^{-2/\alpha}\ell_{1i_{1}}(u)=q_{i_{1}}(u)$ and look at the behavior of
$g(q_{i_{1}}(u))$ when $u\rightarrow\infty.$ We have by choice of
$\mathcal{J}_{0},$
\[
\lim_{u\rightarrow\infty}\frac{q_{i_{1}}(u)}{q_{i}(u)}=\left\{
\begin{array}
[c]{c}%
\kappa_{i},\ \text{if\ }i\in\mathcal{J}_{0};\\
0,\ \text{if }i\notin\mathcal{J}_{0}.
\end{array}
\right.
\]
Hence, by Condition \ref{CondA4},%
\begin{equation}
\lim_{u\rightarrow\infty}u^{2}g(q_{i_{1}}(u))=h(\boldsymbol{\kappa
}),\label{reg1}%
\end{equation}
where $\boldsymbol{\kappa}:=(\kappa_{i},i\in\mathcal{J}_{0};0,i\notin
\mathcal{J}_{0}).$ Now put $st=q_{i_{1}}(u_{1})$; we have, as above,%
\begin{equation}
\lim_{u_{1}\rightarrow\infty}u_{1}^{2}g(q_{i_{1}}(u_{1}))=h(\boldsymbol{\kappa
}).\label{reg2}%
\end{equation}
Using Theorems 1.5.12, 1.5.13 (de Bruijn Lemma) and Proposition 1.5.15,
\cite{BGT}, similarly to \cite{KHP}, we get that for some slowly varying
function $\tilde{\ell}(s)$,
\[
u=s^{-\alpha/2}\tilde{\ell}(s)\ \ \text{and }u_{1}=(st)^{-\alpha/2}\tilde
{\ell}(st),
\]
which together with (\ref{reg1},\ref{reg2}) gives
\[
\lim_{s\rightarrow0}\frac{g(st)}{g(s)}=t^{\alpha},\ \ \alpha=\alpha
(\mathbf{f}),
\]
that is, $r(\mathbf{t})$ regularly varies at direction ${\mathbf{f}}$.
Further, if $\alpha(\mathbf{f})=2$ and $s^{-2}(1-r(s\mathbf{f}))\rightarrow0$
as $s\rightarrow0,$ then by properties of positive defined functions,
$r(t\mathbf{f})\equiv1,$ $t\in\mathbb{R,}$ which contradicts Condition
\ref{CondA3}. Thus Proposition is established.

\begin{remark}
\label{q_for_2} Remark that from Proposition \ref{e_reg_var} it follows that
if $\alpha_{i}=\alpha(\mathbf{e}_{i})=2$ then the corresponding slowly varying
function $\ell_{i}(t)$ is bounded from above. So, since $u^{2}q_{i}^{2}%
(u)\ell_{i}^{2}(q_{i}(u)=\ell_{1i}(u)\ell_{i}^{2}(q_{i}(u)\rightarrow1$ as
$u\rightarrow\infty$, the slowly varied function $\ell_{1i}(u)$ is bounded
from below by a positive constant, that is,%
\begin{equation}
q_{i}(u)\geq q_{0}u^{-1},\ \ q_{0}>0. \label{q for Talagrand}%
\end{equation}
Moreover, it is obvious that (\ref{q for Talagrand}) is valid for any
$q_{i}(u)$, $i=1,...,d$. We shall use this below.
\end{remark}

\begin{remark}
\label{q_i_choice} Observe that for any direction $\mathbf{f}$ we just choose
appropriate $q_{i}(u)$ from the collection of Condition \ref{CondA4}.
\end{remark}

Now assume a behavior of $\sigma(\mathbf{t})$ near its point of absolute
maximum. We shall see from the proof of Lemma \ref{PP-lemma} that the crucial
point is the behavior of the ratio
\[
\frac{1-\sigma^{2}(\mathbf{q}(u)\mathbf{t)}}{1-r(\mathbf{q}(u)\mathbf{t)}}%
\]
as $u\rightarrow\infty.$ In view of Condition \ref{CondA4} we assume the following.

\begin{condition}
\label{CondA5} For any $\mathbf{t}$ there exists the limit
\begin{equation}
h_{1}(\mathbf{t}):=\lim_{u\rightarrow\infty}u^{2}(1-\sigma^{2}(\mathbf{q}%
(u)\mathbf{t))\in\lbrack0,\infty].} \label{cond2}%
\end{equation}

\end{condition}

In case when the limit is equal to zero we speak about the \emph{
stationary-like case}. If the limit is equal to infinity, we refer to the
\emph{Talagrand case}, since M. Talagrand, \cite{talagrand}, has shown that in
most general conditions, for any closed set $S$ and a Gaussian a. s.
continuous function $X(t)$, $t\in S$, having unique point of maximum of
variance, say, at $t_{0}\in S$,
\[
P(S;u)=\mathbf{P}(X(t_{0})>u)(1+o(1)),\ \ u\rightarrow\infty.
\]
In our conditions we show this below. At last, we say about the \emph{
transition case }if $h_{1}(\mathbf{t})$ is neither zero nor infinity. Denote
correspondingly
\begin{align*}
\mathcal{K}_{0} &  :=\{\mathbf{t}\in S\setminus\{\mathbf{0\}:\ }%
h_{1}(\mathbf{t)=}0\},\ \ \mathcal{K}_{c}:=\{\mathbf{t}\in S\setminus
\{\mathbf{0\}:\ }h_{1}(\mathbf{t)\in(}0,\infty)\},\ \ \\
\mathcal{K}_{\infty} &  :=\{\mathbf{t}\in S\setminus\{\mathbf{0\}:\ }%
h_{1}(\mathbf{t)=\infty}\}.
\end{align*}
We shall see that properties of these sets together with all above Conditions
follow asymptotic behavior of the probability $P(S;u)$. Consider one more
simple example which shows that dimensions of $\mathcal{K}_{0},$
$\mathcal{K}_{c},$ $\mathcal{K}_{\infty}$ may be arbitrary. Let $d=2,$
\[
r(t_{1},t_{2})=1-|t_{1}|^{\alpha_{1}}-|t_{2}|^{\alpha_{2}}+o(|\mathbf{t}%
|^{\alpha_{1}}),\mathbf{t\rightarrow0,\ }2\geq\alpha_{1}\geq\alpha_{2}>0,
\]%
\[
\sigma(t_{1},t_{2})=1-|t_{1}|^{\beta_{1}}-|t_{2}|^{\beta_{2}}+o(|\mathbf{t}%
|^{\beta_{1}}),\mathbf{t\rightarrow0,\ }\beta_{1}\geq\beta_{2}>0.
\]
Here one has not to change the basis, having obvious value of $\mathbf{q}(u).$
It is easy to calculate that

\begin{itemize}
\item If $\beta_{1}>\alpha_{1}>\alpha_{2}>\beta_{2}$ then $\dim\mathcal{K}%
_{0}=1,$ $\dim\mathcal{K}_{c}=0,$ $\dim\mathcal{K}_{\infty}=2.$

\item If $\beta_{2}>\alpha_{1}$ then $\dim\mathcal{K}_{0}=2,$ $\dim
\mathcal{K}_{c}=0,$ $\dim\mathcal{K}_{\infty}=0.$

\item If $\beta_{1}=\alpha_{1}>\alpha_{2}>\beta_{2}$ then $\dim\mathcal{K}%
_{0}=0,$ $\dim\mathcal{K}_{c}=1,$ $\dim\mathcal{K}_{\infty}=2.$

\item If $\beta_{1}=\alpha_{1}=\alpha_{2}=\beta_{2}$ then $\dim\mathcal{K}%
_{0}=0,$ $\dim\mathcal{K}_{c}=2,$ $\dim\mathcal{K}_{\infty}=0.$

\item So on.
\end{itemize}

\section{Homogeneous Gaussian fields}

Let $S\subset\mathbb{R}^{d} $ satisfy the assumptions of the previous section,
and $X(\mathbf{t),}$ $\mathbf{t\in}S,$ be an a.s. continuous homogeneous
zero-mean Gaussian field with covariance function satisfying Conditions
\ref{CondA3} and \ref{CondA4}.

\begin{lemma}
\label{lemma_pickands}In the above notations and conditions, for any bounded
closed $T\subset\mathbb{R}^{d}$,%
\[
\mathbf{P}(\max_{\mathbf{t\in q}(u)T}X(\mathbf{t)>}u)=(1+o (1))H_{\mathbf{q}%
}(T\mathbf{)}\Psi(u)
\]
as $u\rightarrow\infty,$ where
\[
H_{\mathbf{q}}(T\mathbf{)=E}\exp(\max_{\mathbf{t\in}T}\chi(\mathbf{t)),\ \ }%
\]
and $\chi(\mathbf{t)}$ is a Gaussian a.s. continuous field with $\chi
(\mathbf{0)=}0,$%
\[
\operatorname*{var}(\chi(\mathbf{t)-}\chi(\mathbf{s))}=2h(\mathbf{t-s),\ E}%
\chi(\mathbf{t)}=-h(\mathbf{t).}%
\]

\end{lemma}

The proof of this lemma is a simple repetition of the proof of Lemma 6.1 of
\cite{book} with applying Condition \ref{CondA4} which implies among other
that $\chi(\mathbf{t})$ exists.

\begin{theorem}
\label{Pickands} In the assumptions of this section, for any closure
$S_{1}\subset S$ of an open set,
\[
\mathbf{P}(\max_{\mathbf{t\in}S_{1}}X(\mathbf{t)>}u)=(1+o(1))|S_{1}%
|H_{\mathbf{q}}\prod_{i=1}^{d}q_{i}^{-1}(u)\Psi(u)
\]
as $u\rightarrow\infty,$ where $|\cdot|$ denotes the volume of $\{\cdot\},$
and
\[
H_{\mathbf{q}}=\lim_{T_{i}\rightarrow\infty,i=1,...,d}(T_{1}...T_{d}%
)^{-1}H_{\mathbf{q}}\left(  \bigotimes_{i=1}^{d}[0,T_{i}]\right)  \in
(0,\infty)\mathbf{.}%
\]
This assertion holds even if $S_{1}$ depends of $u$, $S_{1}=S_{1}(u),$
provided there exist boxes $S_{1}^{\pm}(u)=\otimes_{i=1}^{d}[-S_{1i}^{\pm
}(u),S_{1i}^{\pm}(u)]$ such that $S_{1}^{-}(u)\subset S_{1}(u)\subset
S_{1}^{+}(u)$ with $S_{1i}^{-}(u)q_{i}(u)\rightarrow\infty$ and for some
$\delta<1/2,\ S_{1i}^{+}(u)e^{-\delta u^{2}}\rightarrow0$ as $u\rightarrow
\infty,$ $i=1,...,d.$
\end{theorem}

The proof of the theorem follows step-by-step the proof of Theorem 7.1,
\cite{book}.

\section{Non-homogeneous Gaussian fields}

Now we give two general results for all described above types of behavior of
$\sigma(\mathbf{t})$. The first one is a standard local lemma of Double Sum
Method, a generalization of Lemma \ref{lemma_pickands}, see \cite{book},
\cite{lectures}.

\begin{lemma}
\label{PP-lemma} Under the Conditions $1$ -- $5$, for any $T\in\mathbb{R}%
^{d},$%
\[
P(\mathbf{q}(u)T;u)=(1+o(1))P_{\mathbf{q}}(T\mathbf{)}\Psi(u)
\]
as $u\rightarrow\infty,$ where%
\[
P_{\mathbf{q}}(T\mathbf{)=E}\max_{\mathbf{t}\in T}e^{\chi_{1}(\mathbf{t}%
)},\quad
\]
with $\chi_{1}(\mathbf{t})=\chi(\mathbf{t})-h_{1}(\mathbf{t})$ if
$\mathbf{t}\in\mathcal{K}_{0}\cup\mathcal{K}_{c}$, and $\chi_{1}%
(\mathbf{t})\equiv0$ if $\mathbf{t}\in\mathcal{K}_{\infty}$.
\end{lemma}

The proof of this lemma is a repetition of the corresponding lemma proof in
\cite{book}, using the Conditions \ref{CondA1} -- \ref{CondA5}. The only
essential addition to the proof is careful consideration the points
$\mathbf{t}\in\mathcal{K}_{\infty}$, noticing that the weak convergence of the
field
\[
\chi_{u}(\mathbf{t}):= u(X(\mathbf{q(}u)\mathbf{t})-u)+w
\]
conditioned on $X(\mathbf{0)}=u-w/u$ in $C(T)$ can be restricted to that in
$C(T\setminus\mathcal{K}_{\infty}).$ That is, it can be proved that in the
case of non-empty $\mathcal{K}_{\infty},$
\[
\max_{\mathbf{t}\in T}\chi(\mathbf{t})-h_{1}(\mathbf{t})=\max_{\mathbf{t}\in
T\setminus\mathcal{K}_{\infty}}\chi(\mathbf{t})-h_{1}(\mathbf{t}).
\]

The second general result is extraction of an informative parameter set.
Denote
\begin{equation}
\gamma_{1}(u):=\gamma(u)+u^{-1}\log^{2}u,\label{gamma_1}%
\end{equation}
where $\gamma(u)$ is taken from Proposition \ref{Dmitrovsky}, and
\begin{equation}
B_{u}:=\left\{  \mathbf{t:}1-\sigma^{2}(\mathbf{t}){\leq2}u^{-1}\gamma
_{1}(u)\right\}  .\label{B_u}%
\end{equation}

\begin{lemma}
\label{extracting} In the above conditions and notations,
\[
P(S;u)=\left(  1+O\left(  \exp(-\log^{2}u)\right)  \right)  P(S\cap B_{u};u)
\]
as $u\rightarrow\infty.$
\end{lemma}

\textbf{Proof: }Denote $\sigma_{\varepsilon}:=\sup_{\mathbf{t\in}%
S\setminus\mathbb{B}_{\varepsilon}}\sigma(\mathbf{t}),$ where $\varepsilon$ is
taken from Condition \ref{CondA1}. By Condition \ref{CondA2} we have,
$\sigma_{\varepsilon}<(1+\sigma_{\varepsilon})/2<1.$ By Borell-TIS inequality,
see, for example, \cite{book}, Theorem D.1, for all sufficiently large $u,$
\[
P(S\setminus\mathbb{B}_{\varepsilon};u)\leq2\Psi\left(  \frac{u}%
{(1+\sigma_{\varepsilon})/2}\right)  ,
\]
and the right hand part is exponentially smaller than $P(S;u)$. Indeed,
$\mathbf{0}\in S$ by assumption, hence $P(S;u)\geq\mathbf{P}(X(0)>u)=\Psi
(u)(1+o(1))$ as $u\rightarrow0.$ Further, we have for all sufficiently large
$u,$
\[
\max_{\mathbf{t\in}\mathbb{B}_{\varepsilon}\setminus B_{u}}\sigma
^{2}(\mathbf{t})\leq1-2u^{-1}\gamma_{1}(u),
\]
hence, by Proposition \ref{Dmitrovsky} we have,
\begin{align*}
P(\mathbb{B}_{\varepsilon}\setminus B_{u};u)  &  \leq\exp\left(  -\frac{u^{2}%
}{2(1-2u^{-1}\gamma_{1}(u))}+u\gamma(u)\right)  \leq\exp\left(  -\frac{1}%
{2}u^{2}-u\gamma_{1}(u)+u\gamma(u)\right) \\
&  =\exp\left(  -\frac{1}{2}u^{2}-u\left(  \gamma(u)+u^{-1}\log^{2}u\right)
+u\gamma(u)\right)  =e^{-\log^{2}u}e^{-u^{2}/2}.
\end{align*}
By the above trivial lower estimate, the right hand part is again infinitely
smaller than $P(S;u).$ Thus Lemma is established..

\section{Gaussian fields with unique maximum point of variance.}

We consider first the Gaussian field $X(\mathbf{t)=}\sigma(\mathbf{t)}%
X_{0}(\mathbf{t),}$ satisfying Conditions \ref{CondA1} -- \ref{CondA5}, where
$X_{0}(\mathbf{t})$ is a homogeneous centered Gaussian field with covariance
function satisfying Condition \ref{CondA4} and $\sigma(\mathbf{t)}$ satisfies
Condition \ref{CondA5}$\mathbf{.}$ Then, using standard inequalities including
Slepian Lemma, applied in $B_{u}$, we pass to the general case.

\subsection{Stationary like case.}

Here we consider the stationary-like case, that is $h_{1}(\mathbf{t)=}0$ for
all $\mathbf{t.}$ Fix sufficiently large $u=u_{0}$ and denote
\begin{equation}
f(\mathbf{t})=\frac{1}{2}(1-\sigma^{2}(\mathbf{t})),\ \mathbf{t}\in B_{u_{0}%
}.\label{f(t)}%
\end{equation}
It will be convenient to extend $f(\mathbf{t})$ to the unit cube
$\mathbb{K=\{}\mathbf{t}:\max_{i=1,...,d}|t_{i}|\leq1\},$ having
\[
\max_{\mathbb{K\setminus}B_{u}}f(\mathbf{t)\leq}\max_{\mathbb{\partial
}B_{u_{0}}}f(\mathbf{t})
\]
with $\mathbb{\partial}B,$ the boundary of $B.$ By Lemma \ref{extracting}, the
only behavior of $f(\mathbf{t})$ in any small neighborhood of $\mathbf{0}$
plays role for the desired asymptotic behavior. Hence we assume that
$f(\mathbf{t})$ is continuous in $\mathbb{K}$ and, taking in mind
(\ref{f(t)}),
\begin{equation}
f(\mathbf{t})\in(0,1/2),\ \mathbf{t\in}\mathbb{K\setminus\{}\mathbf{0\}.}%
\label{f(t)1}%
\end{equation}
Introduce the Laplace type integral,%
\begin{equation}
L_{f}(\lambda):=\int_{\mathbb{K}}e^{-\lambda f(\mathbf{t)}}d\mathbf{t}%
,\quad\lambda>0.\label{Lfg}%
\end{equation}
Notice that its asymptotic behavior as $\lambda\rightarrow\infty$ depends only
on behavior $f(\mathbf{t})$ in a vicinity of zero, see, for example,
\cite{fedoruk}.

\begin{proposition}
\label{stationary like} Let Conditions \ref{CondA1}--\ref{CondA5} be
fulfilled. If further $h_{1}(\mathbf{t})=0$ for all $\mathbf{t}$, we have,
\begin{equation}
P(S;u)=(1+o(1))H_{\mathbf{q}}L_{f}(u^{2})\ \prod_{i=1}^{d}q_{i}^{-1}%
(u)\Psi(u),\ \ \label{case S}%
\end{equation}
as $u\rightarrow\infty$.
\end{proposition}

\textbf{Proof}: As it is mentioned above, we consider first a simplified model
for $X,$ that is $X(\mathbf{t})=X_{0}(\mathbf{t})\sigma(\mathbf{t}%
),\ \ \mathbf{t}\in S,$ so that $X(t)$ satisfies Conditions \ref{CondA1} --
\ref{CondA5}$.$ Recall that we consider the case
\begin{equation}
\lim_{u\rightarrow\infty}\frac{1-\sigma^{2}(\mathbf{q(}u\mathbf{)t}%
)}{1-r(\mathbf{q(}u\mathbf{)t})}=0.\label{case S1}%
\end{equation}
Denote
\begin{equation}
T_{i}(u)=\sup\{|t_{i}|:\mathbf{t}\in B_{u}\},\ \text{ }i=1,...,d.\label{T+T-}%
\end{equation}
Obviously that for all $i$, $T_{i}(u)>0$ and it tends to zero as
$u\rightarrow\infty.$ Moreover,
\[
\bigotimes_{i=1}^{d}[-T_{i}(u),T_{i}(u)]\supseteq B_{u}.
\]
By (\ref{case S1}), (\ref{T+T-}) and the definition of $B_{u}$,
\[
\lim_{u\rightarrow\infty}\frac{T_{i}(u)}{q_{i}(u)}=\infty,\ \ i=1,...,d.
\]
Let an increasing $\kappa(u)$ be such that
\begin{equation}
\lim_{u\rightarrow\infty}\kappa(u)=\lim_{u\rightarrow\infty}\frac{T_{i}%
(u)}{q_{i}(u)\kappa(u)}=\infty,\ \ i=1,...,d.\label{kappa}%
\end{equation}
The box
\begin{equation}
\Delta_{\mathbf{0}}:=\kappa(u)\bigotimes_{i=1}^{d}[0,q_{i}(u)]\label{Delta_0}%
\end{equation}
satisfies conditions Theorem \ref{Pickands}, hence,
\begin{equation}
\mathbf{P}\Bigl(\max_{\mathbf{t}\in\Delta_{\mathbf{0}}}X_{0}(\mathbf{t}%
)>u\Bigr)=(1+\delta(u))H_{\mathbf{q}}|\Delta_{\mathbf{0}}|\prod_{i=1}^{d}%
q_{i}^{-1}(u)\Psi(u),\label{P(A0)}%
\end{equation}
with $\delta(u)\rightarrow0$ as $u\rightarrow\infty.$ Denote
\[
\Delta_{\mathbf{k}}(u)=\kappa(u)\mathbf{kq}(u)+\Delta_{\mathbf{0}%
}(u)\mathbf{,}\quad\mathbf{k}\in\mathbb{Z}^{d},u>0.
\]
For all $\mathbf{k}$ with $\Delta_{\mathbf{k}}(u)\cap B_{u}\neq\varnothing$
introduce events
\[
A_{\mathbf{k}}(u)=\Bigl\{\max_{\mathbf{t}\in\Delta_{\mathbf{k}}(u)}%
X_{0}(t)>u_{\mathbf{k}}\Bigr\}\ \text{with }u_{\mathbf{k}}=u/\sigma
_{\mathbf{k}},\ \sigma_{\mathbf{k}}=\max_{\mathbf{t}\in\Delta_{\mathbf{k}}%
(u)}\sigma(\mathbf{t}),
\]
and
\[
A_{\mathbf{k}}^{\prime}(u)=\Bigl\{\max_{\mathbf{t}\in\Delta_{\mathbf{k}}%
(u)}X_{0}(\mathbf{t})>u_{k}^{\prime}\Bigr\}\ \text{with }u_{\mathbf{k}%
}^{\prime}=u/\sigma_{\mathbf{k}}^{\prime},\ \sigma_{\mathbf{k}}^{\prime}%
=\min_{\mathbf{t}\in\Delta_{\mathbf{k}}(u)}\sigma(\mathbf{t}).
\]
By definition of $B_{u},$ after some easy calculations we have that
\begin{equation}
u\leq u_{\mathbf{k}},u_{\mathbf{k}}^{\prime}\leq u+\gamma_{1}(u),\ \mathbf{k}%
\in K_{u}:=\{\mathbf{k}:\Delta_{\mathbf{k}}(u)\cap B_{u}\neq\varnothing
\}.\label{K_u}%
\end{equation}
Hence all the boxes $\Delta_{\emph{k}}(u)$ satisfy Theorem \ref{Pickands}
conditions with $u_{\mathbf{k}}$ instead of $u.$ Therefore
\[
\mathbf{P}(A_{\mathbf{k}}(u))=(1+\delta(u_{\mathbf{k}}))H_{\mathbf{q}}%
|\Delta_{\mathbf{k}}|\prod_{i=1}^{d}q_{i}^{-1}(u_{\mathbf{k}})\Psi
(u_{\mathbf{k}}),
\]
and
\[
\mathbf{P}(A_{\mathbf{k}}^{\prime}(u))=(1+\delta(u_{\mathbf{k}}^{\prime
}))H_{\mathbf{q}}|\Delta_{\mathbf{k}}|\prod_{i=1}^{d}q_{i}^{-1}(u_{\mathbf{k}%
}^{\prime})\Psi(u_{\mathbf{k}}^{\prime}),
\]
$\ \mathbf{k}\in K_{u}.$ By definition (\ref{K_u}) of $K_{u},$ there exists a
positive non-increasing $\delta_{1}(u)$ tending to zero as $u\rightarrow
\infty,$ such that for all $\mathbf{k}\in K_{u},$
\[
1-\delta_{1}(u)\leq\frac{\min(q_{i}(u_{\mathbf{k}}),q_{i}(u_{\mathbf{k}%
}^{\prime}))}{q_{i}(u)}\leq\frac{\max(q_{i}(u_{\mathbf{k}}),q_{i}%
(u_{\mathbf{k}}^{\prime}))}{q_{i}(u)}\leq1+\delta_{1}(u),\ \ i=1,...,d.
\]
Further, since $\delta_{2}(u)=\sup_{v\geq u}|\delta(v)|\rightarrow0$ as
$u\rightarrow\infty$ and $u_{\mathbf{k}},u_{\mathbf{k}}^{\prime}\geq u,$ we
get that for all $\mathbf{k}\in K_{u},$
\begin{equation}
\mathbf{P}(A_{\mathbf{k}}(u)),\mathbf{P}(A_{\mathbf{k}}^{\prime}%
(u))\lesseqgtr(1\pm\delta_{1}(u))(1\pm\delta_{2}(u))H_{\mathbf{q}}%
|\Delta_{\mathbf{k}}|\prod_{i=1}^{d}q_{i}^{-1}(u)\Psi(u_{\mathbf{k}%
}).\label{P(A_ku)}%
\end{equation}
By Bonferroni inequality%
\begin{equation}
{P(B}_{u};u)\leq\sum_{\mathbf{k}:\Delta_{\mathbf{k}}(u)\cap B_{u}%
\neq\varnothing}\mathbf{P}(A_{k}(u)),\label{Bonf1}%
\end{equation}
and%
\begin{equation}
{P(B}_{u};u)\geq\sum_{\mathbf{k}:\Delta_{\mathbf{k}}(u)\subset B_{u}%
}\mathbf{P}(A_{\mathbf{k}}^{\prime}(u))-\sum_{\mathbf{k,l}:\Delta_{\mathbf{k}%
}(u)\cap B_{u}\neq\varnothing,\Delta_{\mathbf{l}}(u)\cap B_{u}\neq
\varnothing,\mathbf{k\neq l}}\mathbf{P}(A_{\mathbf{k}}(u)A_{\mathbf{l}%
}(u)).\label{Bonf2}%
\end{equation}
Write
\begin{equation}
u_{\mathbf{k}}^{2}=u^{2}+u^{2}(1-\sigma_{\mathbf{k}}^{2})+\frac{u^{2}%
(1-\sigma_{\mathbf{k}}^{2})^{2}}{\sigma_{\mathbf{k}}^{2}}.\label{s1}%
\end{equation}
Using
\begin{equation}
\frac{u^{2}(1-\sigma_{\mathbf{k}}^{2})^{2}}{2\sigma_{\mathbf{k}}^{2}}\geq
\frac{u^{2}u^{-2}\gamma_{1}^{2}(u)}{2},\label{s2}%
\end{equation}
since $\gamma_{1}(u)\rightarrow0$ as $u\rightarrow\infty,$ we have for some
positive $\delta_{3}(u)$ with $\delta_{3}(u)\rightarrow0,$ $u\rightarrow
\infty,$ uniformly in $\mathbf{k}\in K_{u},$%
\begin{equation}
\sigma_{\mathbf{k}}\exp\left(  -\frac{u^{2}(1-\sigma_{\mathbf{k}}^{2})^{2}%
}{2\sigma_{\mathbf{k}}^{2}}\right)  \lesseqgtr1\pm\delta_{3}(u).\label{s3}%
\end{equation}
Hence, by (\ref{P(A_ku)})
\begin{equation}
\sum_{\mathbf{k}:\Delta_{\mathbf{k}}(u)\cap B_{u}\neq\varnothing}%
\mathbf{P}(A_{k}(u))\lesseqgtr(1\pm\delta_{4}(u))H_{\mathbf{q}}\prod_{i=1}%
^{d}q_{i}^{-1}(u)\Psi(u)\sum_{\mathbf{k}:\ \Delta_{\mathbf{k}}\cap B_{u}%
\neq\varnothing}|\Delta_{\mathbf{k}}|e^{-u^{2}(1-\sigma_{\mathbf{k}}^{2}%
)/2},\label{s4}%
\end{equation}
where
\[
1\pm\delta_{4}(u)=(1\pm\delta_{1}(u))(1\pm\delta_{2}(u))(1\pm\delta_{3}(u)).
\]
Remark that all the relations (\ref{s1} -- \ref{s3}) are also valid for for
$\sigma_{\mathbf{k}}^{\prime},$ with some other $\delta_{\nu}(u),$
$\nu=1,2,3,4,$ having the same properties, say, for $\delta_{\nu}^{\prime
}(u),$ $\nu=1,2,3,4.$ Hence we also have (\ref{s4}) with $A_{k}^{\prime}(u)$
and $\delta_{4}^{\prime}(u)$ instead of $A_{k}(u)$ and $\delta_{4}(u).$
Denote
\begin{equation}
\Sigma(u):=\sum_{\mathbf{k}:1-\sigma_{\mathbf{k}}^{2}\leq2u^{-1}\gamma_{1}%
(u)}|\Delta_{\mathbf{k}}|e^{-u^{2}(1-\sigma_{\mathbf{k}}^{2})/2}%
,\ \ \text{and\ }\Sigma^{\prime}(u):=\sum_{\mathbf{k}:1-\sigma_{\mathbf{k}%
}^{\prime2}\leq2u^{-1}\gamma_{1}(u)}|\Delta_{\mathbf{k}}|e^{-u^{2}%
(1-\sigma_{\mathbf{k}}^{\prime2})/2}.\label{s5}%
\end{equation}
Notice that $|\Delta_{\mathbf{k}}|=|\Delta_{\mathbf{0}}|=\kappa^{d}%
(u)\prod_{i=1}^{d}q_{i}(u),$ we used $|\Delta_{\mathbf{k}}|$ above for
visibility. We have that $\Sigma(u)$ and $\Sigma^{\prime}(u)$ are integral
sums for the integral%
\[
I(u):=\int_{f(\mathbf{t})\leq2u^{-1}\gamma_{1}(u)}e^{-u^{2}f(\mathbf{t}%
)}d\mathbf{t}.
\]
Besides,
\[
\Sigma^{\prime}(u)\leq I(u)\leq\Sigma(u),
\]
therefore from inequalities (\ref{s4}) and mentioned there their counterparts
for $A_{k}^{\prime}(u)$ and $\delta_{4}^{\prime}(u)$ it follows that for some
$\tilde{\delta}(u)$ tending to zero as $u\rightarrow\infty,$%
\[
\sum_{\mathbf{k}:\Delta_{\mathbf{k}}(u)\cap B_{u}\neq\varnothing}%
\mathbf{P}(A_{k}(u))\lesseqgtr(1\pm\tilde{\delta}(u))H_{\mathbf{q}}\prod
_{i=1}^{d}q_{i}^{-1}(u)\Psi(u_{\mathbf{k}})\int_{f(\mathbf{t})\leq
2u^{-1}\gamma_{1}(u)}e^{-u^{2}f(\mathbf{t})}d\mathbf{t}.
\]
Now, using given by (\ref{f(t)} -- \ref{Lfg}) definitions of $f(\mathbf{t})$
and $L(\lambda),$ we have,
\[
\int_{f(\mathbf{t})\leq2u^{-1}\gamma_{1}(u)}e^{-u^{2}f(\mathbf{t})}%
d\mathbf{t}=\left(  1+O\left(  e^{-2\log^{2}u}\right)  \right)  L_{f}%
(u^{2})\mathbf{.}%
\]
Indeed, it follows from (\ref{gamma_1}) that $\gamma_{1}(u)\geq u^{-1}\log
^{2}u,$ so that for all sufficiently small $\mathbf{t}$ lying outside of the
integration domain, $u^{2}f(\mathbf{t})\geq2u^{2}u^{-1}\gamma_{1}(u)\geq
2u^{2}u^{-2}\log^{2}u.$

Now estimate from above the double sum in (\ref{Bonf2}). We have for non
neighboring boxes,
\begin{align}
\mathbf{P}(A_{\mathbf{k}}A_{\mathbf{l}}) &  \leq\mathbf{P}(\max_{(\mathbf{s,t}%
)\in\Delta_{\mathbf{k}}\otimes\Delta_{\mathbf{l}}}X(\mathbf{s})+X(\mathbf{t}%
)>u_{\mathbf{k}}+u_{\mathbf{l}})\nonumber\\
&  \leq\mathbf{P}\left(  \max_{(\mathbf{s,t})\in\Delta_{\mathbf{k}}%
\otimes\Delta_{\mathbf{l}}}\frac{X_{0}(\mathbf{s})+X_{0}(\mathbf{t})}%
{\sqrt{2+2r(\mathbf{t-s})}}\geq\frac{u_{\mathbf{k}}+u_{\mathbf{l}}}%
{\sqrt{2(1+r_{\mathbf{k,l}}(u)}}\right)  ,\label{DS1}%
\end{align}
where $2+2r(\mathbf{t-s})$ is the variance of $X_{0}(\mathbf{s})+X_{0}%
(\mathbf{t})$ and
\[
r_{\mathbf{k,l}}(u)=\max_{(\mathbf{s,t})\in\Delta_{\mathbf{k}}\otimes
\Delta_{\mathbf{l}}}r(\mathbf{t-s}).
\]
For the increments of the zero mean field
\[
Y(\mathbf{s,t}):=\frac{X_{0}(\mathbf{s})+X_{0}(\mathbf{t})}{\sqrt
{2+2r(\mathbf{t-s})}}%
\]
with unit variance, one can get by simple algebra that for an absolute
constant $C,$
\[
\mathbf{E}(Y(\mathbf{s,t})-Y(\mathbf{s}^{\prime}\mathbf{,t}^{\prime})^{2}\leq
C(\mathbf{E}(X_{0}(\mathbf{s})-X_{0}(\mathbf{s}^{\prime}))^{2}+\mathbf{E}%
(X_{0}(\mathbf{t})-X_{0}(\mathbf{t}^{\prime}))^{2}).
\]
From this inequality by standard Gaussian technique including Slepian
inequality it follows that the right hand part of (\ref{DS1}) is at most
\[
C_{1}\kappa^{2d}(u)\Psi\left(  \frac{u_{\mathbf{k}}+u_{\mathbf{l}}}%
{\sqrt{2(1+r_{\mathbf{k,l}}(u)}}\right)  ,
\]
where the constant $C_{1}$ does not depend of $\mathbf{k,l}.$ Now write,
\[
(u_{\mathbf{k}}+u_{\mathbf{l}})^{2}=u^{2}\left(  \sigma_{\mathbf{k}}%
^{-1}+\sigma_{\mathbf{l}}^{-1}\right)  ^{2}=u^{2}\left(  2+(1-\sigma
_{\mathbf{k}})+(1-\sigma_{\mathbf{l}})+\sum_{j=2}^{\infty}((1-\sigma
_{\mathbf{k}})^{j}+(1-\sigma_{\mathbf{l}})^{j})\right)  ^{2}.
\]
Using that $r(\mathbf{t})$ regularly varies at zero in any direction
$\mathbf{f}$, $|\mathbf{f}|=1$, with indexes $\alpha(\mathbf{f})\in(0,2]$, see
(\ref{alpha}), we get for some positive $\gamma\in\lbrack\max_{\mathbf{f}%
}\alpha(\mathbf{f}),2]$ and $\Gamma$ that $1-r(\mathbf{t})\geq\Gamma
|\mathbf{t}|^{\gamma}.$ Remark that from (\ref{A3}) and followed then argument
in case $\alpha(\mathbf{f})=2$ the corresponding slowly varying function is
bounded at zero otherwise square mean derivative of $X_{\mathbf{0}}$ in this
direction exists and is a constant, this contradicts the conditions on $r$ and
$R.$ Thus we have for non neighboring $\mathbf{k,l}$,
\[
u^{2}(1-r_{\mathbf{k,l}}(u))\geq\Gamma\kappa^{\gamma}(u)u^{2-2/\gamma}%
\geq\Gamma\kappa^{\gamma}(u).
\]
Therefore, by analogy to (\ref{s1}, \ref{s2}), we have,%
\begin{align*}
\frac{(u_{\mathbf{k}}+u_{\mathbf{l}})^{2}}{2(1+r_{\mathbf{k,l}}(u))} &
\geq\frac{4u^{2}}{2(2-(1-r_{\mathbf{k,l}}(u))}+\frac{1}{2}u^{2}(1-\sigma
_{\mathbf{k}})+\frac{1}{2}u^{2}(1-\sigma_{\mathbf{l}})\\
&  \geq u^{2}+\frac{1}{2}\Gamma\kappa^{\gamma}(u)+\frac{1}{4}u^{2}%
(1-\sigma_{\mathbf{k}}^{2})+\frac{1}{4}u^{2}(1-\sigma_{\mathbf{l}}^{2}),
\end{align*}
where we also used that $r\leq1$ and $\sigma\leq1.$ Finally we have for non
neighboring $\mathbf{k,l},$ that is, $\max_{i=1,...,d}|k_{i}-l_{i}|>1,$%
\[
\mathbf{P}(A_{\mathbf{k}}A_{\mathbf{l}})\leq C_{2}\kappa^{2}(u)\exp\left(
-\frac{1}{2}\Gamma|\kappa(u)|^{\gamma}\right)  \Psi(u)\exp\left(  -\frac{1}%
{4}u^{2}(1-\sigma_{\mathbf{k}}^{2})-\frac{1}{4}u^{2}(1-\sigma_{\mathbf{l}}%
^{2})\right)  .
\]
Summing up, denoting $\mathcal{A(}\mathbf{k,l}):=\{\mathbf{k,l:}%
\Delta_{\mathbf{k}}(u)\cap B_{u}\neq\varnothing,\Delta_{\mathbf{l}}(u)\cap
B_{u}\neq\varnothing\}$, we have for some $C_{2}>0$,
\begin{align}
\sum_{\mathbf{k,l}\in\mathcal{A(}\mathbf{k,l}):\max_{i=1,...,d}|k_{i}%
-l_{i}|>1} &  \mathbf{P}(A_{\mathbf{k}}A_{\mathbf{l}})\leq C_{2}\kappa
^{2}(u)\exp(-\Gamma|\kappa(u)|^{\gamma})\Psi(u)\nonumber\\
&  \times\left(  \sum_{\mathbf{k}:\Delta_{\mathbf{k}}(u)\cap B_{u}%
\neq\varnothing}e^{-\frac{1}{4}u^{2}(1-\sigma_{\mathbf{k}}^{2})}\right)
^{2}.\label{DS2}%
\end{align}
By the above argument the last sum multiplied by $\kappa^{d}(u)\prod_{i=1}%
^{d}q_{i}(u)$ is an integral sum for the integral $I(u/\sqrt{2}),$ with
inessential changing of the integration domain. Obvious application of Schwarz
inequality gives
\[
I^{2}(u/\sqrt{2})\leq(1+\varepsilon)I(u),
\]
for any $\varepsilon>0$ and all sufficiently large $u.$ The first exponent in
the right hand part of (\ref{DS2}) gives that the double sum over non
neighboring intervals is infinitely smaller than both the single sums in
(\ref{Bonf1}, \ref{Bonf2}).

Consider the double sum over neighboring intervals, that is over
$\mathbf{k,l}$ with $\max_{i=1,...,d}|k_{i}-l_{i}|=1.$ This part is quite
similar to the corresponding argument in \cite{book}, \cite{lectures}. Let for
definiteness $l_{1}-k_{1}=1.$ Denote
\[
\Delta_{\mathbf{0}}^{\prime}:=\left[  0,\sqrt{\kappa(u)q_{1}(u)}\right]
\otimes\bigotimes_{i=2}^{d}[0,\kappa(u)q_{i}(u)]
\]
and write
\begin{align*}
\mathbf{P}(A_{\mathbf{k}}A_{\mathbf{l}})  &  \leq\mathbf{P}(\max
_{\mathbf{t}\in\kappa(u)\mathbf{lq(u)+}\Delta_{\mathbf{0}}^{\prime}}%
X_{0}(\mathbf{t})\geq u_{\mathbf{l}})\\
&  +\mathbf{P}(\max_{\mathbf{s}\in\Delta_{\mathbf{k}}\mathbf{,t}%
\in\mathbf{lq(u)}+\Delta_{\mathbf{0}}\setminus\Delta_{\mathbf{0}}^{\prime}%
}X_{0}(\mathbf{s})+X_{0}(\mathbf{t})\geq u_{\mathbf{k}}+u_{\mathbf{l}}).
\end{align*}
The sum of the first probabilities on the right can be estimated using the
same argument as the estimation of single sums above. Wherein the multiplier
$\sqrt{\kappa(u)q_{1}(u)}$ appears which gives that the sum is infinitely
smaller the single sum above. For the second probability on the right the
argument of the double sum estimation over non-neighboring boxes can be
applied because of the distance between boxes is not zero but $\sqrt
{\kappa(u)q_{1}(u)}.$ This also gives that the sum is infinitely smaller than
the single sums.

Hence in view of already mentioned standard passage from the particular
$X(\mathbf{t})=X_{0}(\mathbf{t})\sigma(\mathbf{t})$ to the general Gaussian
process (by applying Slepian inequality), the proof follows. Thus Proposition
is established.

\begin{remark}
\label{DS_r} In the case when the fraction (\ref{case S1}) tends to zero
sufficiently fast, for example for $d=1,$
\[
\lim_{t\rightarrow0}\frac{1-\sigma^{2}(t)}{t^{\varepsilon}(1-r(t))}=0
\]
with some $\varepsilon>0,$ the estimation of the double sum is quite similar
to that in \cite{book}, \cite{lectures}. But the fraction may tend to zero
very slowly, so that for this situation the evaluations and estimations have
to be more precise, what we have done here.
\end{remark}

\subsection{Talagrand case.}

\begin{proposition}
\label{Talagrand_case}In the above conditions, if $\mathcal{K}_{\infty
}=S\setminus\{\mathbf{0}\},$ then
\begin{equation}
P(S;u)=(1+o(1))\Psi(u),\ \ \ \label{case T}%
\end{equation}
as $u\rightarrow\infty.$
\end{proposition}

\textbf{Proof:} By Proposition \ref{e_reg_var} and Remark \ref{q_for_2}, in
view of Lemma \ref{extracting}, we have that for $\varepsilon\in(0,q_{0}),$
and all sufficiently large $u$, $B_{u}\subseteq\varepsilon\bigotimes_{i=1}%
^{d}[-q_{i}(u),q_{i}(u)].$ Further, similarly to the proof of Lemma 8.4,
\cite{book}, we get that
\[
\limsup_{u\rightarrow\infty}\frac{P(B_{u};u)}{\Psi(u)}\leq\limsup
_{u\rightarrow\infty}\frac{P(\varepsilon\bigotimes_{i=1}^{d}[-q_{i}%
(u),q_{i}(u)];u)}{\Psi(u)}\leq\mathbf{E}\max_{\mathbf{t}\in\lbrack
-\varepsilon,\varepsilon]^{d}}e^{\chi(\mathbf{t})}.
\]
Observe that in the latter inequality we again pass to a homogeneous field
using monotonicity with respect to the variance and Slepian's inequality, with
following application of Lemma \ref{lemma_pickands}. Then, as in the proof of
Lemma 8.4, we use Monotone Convergence Theorem to let $\varepsilon
\downarrow0.$

\subsection{The transition case}

Asymptotic evaluations for the exceeding probabilities in this case are quite
similar to the corresponding evaluations in \cite{PP}, \cite{book},
\cite{lectures}, with applying Lemma \ref{PP-lemma}. Take in the Lemma
\[
T=T_{\mathbf{0}}=[-K,K]^{d},
\]
and denote
\[
T_{\mathbf{k}}=2\mathbf{k}K+T_{\mathbf{0}},\ \mathbf{k\in}\mathbb{Z}^{d}.
\]
First we apply Lemma \ref{PP-lemma} for $T_{\mathbf{0}},$ then we estimate
from above the sum of the probabilities $P(\mathbf{q}(u)T_{\mathbf{k}}\cap
B_{u};u)$ over $\mathbf{k\neq0}$ with $\mathbf{q}(u)T_{\mathbf{k}}\cap
B_{u}\neq\varnothing$ using the same Lemma and regular varying of
$r(\mathbf{t})$ in any direction. Then we let $K$ tend to infinity. On this
way we get the following.

\begin{proposition}
\label{PPcase}In the above conditions, if $\mathcal{K}_{c}=S\setminus
\{\mathbf{0}\},$ then%
\begin{equation}
P(S;u)=(1+o(1))P_{\mathbf{q}}\Psi(u) \label{case P}%
\end{equation}
as $u\rightarrow\infty,$ with $P_{\mathbf{q}}=\lim_{K\rightarrow\infty
}P_{\mathbf{q}}([-K,K]^{d})\in(0,\infty).$
\end{proposition}

\subsection{General case.}

First formulate several simple generalizations of above propositions. The
first one is a generalization of Proposition \ref{Talagrand_case}.

\begin{proposition}
\label{K_infty} In the above conditions,
\[
P(S;u)=P(S\setminus\mathcal{K}_{\infty};u)(1+o(1))\text{ as }u\rightarrow
\infty.
\]

\end{proposition}

The proof repeats the proof of Proposition \ref{Talagrand_case}, with Lemma
\ref{extracting} application. The second one is a simple reformulation of
Proposition \ref{PPcase}.

\begin{proposition}
\label{PP_gen} In the above conditions, if $\dim\mathcal{K}_{c}>0,$ then
\begin{equation}
P(\mathcal{K}_{c};u)=(1+o(1))P_{\mathbf{q}}(\mathcal{K}_{c})\Psi(u),
\label{K_c_fin}%
\end{equation}
as $u\rightarrow\infty,$ where
\[
P_{\mathbf{q}}(\mathcal{K}_{c})=\lim_{K\rightarrow\infty}\mathbf{E}%
\max_{\mathbf{t}\in\lbrack-K,K]^{d}\cap\mathcal{K}_{c}}e^{\chi_{1}%
(\mathbf{t})}\in(0,\infty).
\]

\end{proposition}

The proof starts with the set%
\[
T=T_{\mathbf{0}}=[-K,K]^{d}\cap\mathcal{K}_{c},
\]
with followed corresponding definition of $T_{\mathbf{k}}.$

The proof of next proposition repeats the proof of Proposition
\ref{stationary like}.

\begin{proposition}
\label{dimK_0=d}If $\dim\mathcal{K}_{0}=d$, then
\[
P(\mathcal{K}_{0};u)=(1+o(1))H_{\mathbf{q}}L_{f}(\mathcal{K}_{0},u^{2}%
)\ \prod_{i=1}^{d}q_{i}^{-1}(u)\Psi(u),\ \
\]
as $u\rightarrow\infty,$ where
\begin{equation}
L_{f}(\mathcal{K}_{0},\lambda):=\int_{\mathcal{K}_{0}\cap\mathbb{K}%
}e^{-\lambda f(\mathbf{t)}}d\mathbf{t},\quad\lambda>0,\label{L_gen}%
\end{equation}
and $f(\mathbf{t})$ is given by (\ref{f(t)}, \ref{f(t)1}).
\end{proposition}

Recall that $\mathbb{K}$ is the unit cube, see (\ref{Lfg}).

Notice that the case $\dim\mathcal{K}_{0}=0$ means $\mathcal{K}_{0}%
=\emptyset,$ by definition. Assume now that $\dim\mathcal{K}_{0}>0.$
Generally, since the behavior of $\sigma(\mathbf{t})$ near its maximum point
$\mathbf{0}$ can be various, $\mathcal{K}_{0}$ may consist of several
non-intersecting connected manifolds of various dimensions. In order to avoid
technical difficulties due to too exotic behavior of $\sigma(\mathbf{t),}$
assume the following.

\begin{condition}
\label{CondA6} The set $\mathcal{K}_{0}$ consists of finite number of smooth
(two times continuously differentiable) disjoint manifolds, namely,
\begin{equation}
\mathcal{K}_{0}=\bigcup_{i=1}^{n}\mathcal{K}_{0i},\ \dim\mathcal{K}_{0i}%
=k_{i},\ 0<k_{1}\leq k_{2}\leq...\leq k_{n}\leq d. \label{part}%
\end{equation}
Assume that for any $i,$ $k_{i}$-dimensional volume of $\mathcal{K}_{0i}$ is
finite, $|\mathcal{K}_{0i}|<\infty.$
\end{condition}

Fix $i$ with $k_{i}<d,$ and consider in $\mathcal{K}_{0i}$ curvilinear
coordinates. For $\mathbf{t\in}\mathcal{K}_{0i},$ using Proposition
\ref{e_reg_var}, choose coordinate vectors $\mathbf{e}_{j}(\mathbf{t),}$
$j=1,...,k_{i}\mathbf{\ }$of this curvilinear coordinates and complete them to
a basis  $\{\mathbf{e}_{j}(\mathbf{t),}$ $j=1,...,k_{i},$ $\mathbf{\tilde{e}%
}_{j}(\mathbf{t),}$ $j=k_{i}+1,...,d\}$ in $\mathbb{R}^{d}.$ and denote
\begin{equation}
\mathbf{q}^{i,\mathbf{t}}(u):=(q_{j}^{i,\mathbf{t}}(u),\ j=1,...,k_{i}%
,\tilde{q}_{j}^{i,\mathbf{t}}(u),j=k_{i}+1,...,d),\label{K01}%
\end{equation}
with corresponding positive limits
\begin{equation}
\lim_{u\rightarrow\infty}u^{2}(1-r(\mathbf{q}^{i,\mathbf{t}}(u)\mathbf{s}%
)=:h^{i,\mathbf{t}}(\mathbf{s}),\label{K02}%
\end{equation}
where $\mathbf{s}=(s_{1},...,s_{d})$ is written in these coordinates. Remark
that by Proposition \ref{e_reg_var}, functions $q_{j}^{i,\mathbf{t}}(u)$ are
taken from the collection of Condition \ref{CondA4}, but the choice of them
can depend on $\mathbf{t}$ and $i,$ the index of the manifold. In fact,
$h^{i,\mathbf{t}}(\mathbf{s})=h(U_{i,\mathbf{t}}\mathbf{s),}$ where
$U_{i,\mathbf{t}}$ is an orthogonal transition matrix to the curvilinear
coordinates with the above  orthogonal complement to a basis in $\mathbb{R}%
^{d}.$

We have,%
\begin{equation}
\lim_{u\rightarrow\infty}\frac{1-\sigma(\mathbf{q}^{i,\mathbf{t}}%
(u)\mathbf{s})}{1-r(\mathbf{q}^{i,\mathbf{t}}(u)\mathbf{s})}=0. \label{K03}%
\end{equation}
By analogy with relations (\ref{case S1} - \ref{kappa}) in Proposition
\ref{stationary like} proof, we build a partition of $\mathcal{K}_{0i}$ with
similar to $\Delta_{0}$ $k_{i}$-dimensional blocks, denote them by
$\Delta(\mathbf{t}_{\nu}),$ where $\{\mathbf{t}_{\nu},\nu=1,...,N\}$ is a grid
satisfying (\ref{K01} - \ref{K03}) for all $\mathbf{s\in}\Delta(\mathbf{t}%
_{\nu}),$ $\nu=1,...,N,$ correspondingly, that is,
\begin{equation}
\mathcal{K}_{0i}=\bigcup_{i=1}^{N}\Delta(\mathbf{t}_{\nu}). \label{grid}%
\end{equation}
Similarly to the proof of Proposition \ref{stationary like}, but for
$\mathcal{K}_{0i}\cap B_{u}$ instead of $S\cap B_{u},$ using Theorem
\ref{Pickands} to get (\ref{P(A0)}) for all $\Delta(\mathbf{t}_{\nu})$ and
thicken the grid unboundedly, we get, using Condition \ref{CondA6}, the
following Lemma.

\begin{lemma}
\label{K0i} For any $\mathcal{K}_{0i}$ from the partition (\ref{part}) with
$k_{i}<d,$
\[
P(\mathcal{K}_{0i};u)=(1+o(1))\int_{\mathcal{K}_{0i}}H_{\mathbf{q}%
^{i,\mathbf{t}}}\prod_{i=1}^{k_{i}}(q_{j}^{i,\mathbf{t}}(u)^{-1}%
e^{-u^{2}f(\mathbf{t)}}\nu_{k_{i}}(d\mathbf{t)})\Psi(u)\ \
\]
as $u\rightarrow\infty$, where $\nu_{k_{i}}(d\mathbf{t)}$ is an elementary
$k_{i}$-dimensional volume of $\mathcal{K}_{0i},$ and $f(\mathbf{t})$ is given
by (\ref{f(t)}).
\end{lemma}

\begin{remark}
Notice that if a manifold $\mathcal{K}_{0i}$ is a linear subspace, all
$\mathbf{q}^{i,\mathbf{t}}(u)$ do not depend of $\mathbf{t,}$ so that%
\[
P(\mathcal{K}_{0i};u)=(1+o(1))H_{\mathbf{q}^{i}}\prod_{i=1}^{k_{i}}(q_{j}%
^{i}(u)^{-1}\int_{\mathcal{K}_{0i}}e^{-u^{2}f(\mathbf{t)}}\nu_{k_{i}%
}(d\mathbf{t)})\Psi(u)\ \
\]
as $u\rightarrow\infty$.
\end{remark}

Turning to (\ref{part}), we have,

\begin{proposition}
\label{Th_dimK_0<d} If $\dim\mathcal{K}_{0}<d,$
\begin{equation}
P(\mathcal{K}_{0};u)=(1+o(1))\sum_{i=1}^{n}\int_{\mathcal{K}_{0i}%
}H_{\mathbf{q}^{i,\mathbf{t}}}\prod_{i=1}^{k_{i}}(q_{j}^{i,\mathbf{t}}%
(u)^{-1}e^{-u^{2}f(\mathbf{t)}}\nu_{k_{i}}(d\mathbf{t)})\Psi
(u)\ \ \label{K0_fin2}%
\end{equation}
as $u\rightarrow\infty$.
\end{proposition}

Remark that the summands in (\ref{K0_fin2}) can have different orders in $u$
depending on the dimension of the corresponding component $\mathcal{K}_{0i},$
on behavior of $q_{j}^{i,\mathbf{t}}(u)$s and on behavior of $\sigma
^{2}(\mathbf{t})$. Hence only summands with slowest order play a role. Remark
also that by Proposition \ref{dimK_0=d}, if for some $i,$ $\dim\mathcal{K}%
_{0i}=d,$ no summands in (\ref{K0_fin2}) contribute to the asymptotics of
$P(S;u).$

\textbf{Proof:} We have only to estimate the double probabilities
$\mathbf{P}(A(\Delta(\mathbf{s}_{\nu}))A(\Delta(\mathbf{t}_{\mu}))),$ where
$A(\Delta(\mathbf{s}_{\nu}))$ and $A(\Delta(\mathbf{t}_{\mu}))$ are events
generated by corresponding partitions in different component manifolds. In
view of Proposition \ref{e_reg_var}, denote
\[
\mathbb{K}_{\alpha/2}:=u^{-4/\alpha}\mathbb{K=\{}\mathbf{t:}|t_{i}|\leq
u^{-4/\alpha}\},
\]
where $\alpha$ is defined in (\ref{alpha}), the minimal index of regular
variation of $r(t\mathbf{f})$ at zero, $|\mathbf{f|}=1.$ Then for the sets
$\Delta(\mathbf{t})$ from partitions of the component manifolds containing in
the manifold $\mathcal{K}_{0}\setminus\mathbb{K}_{\alpha/2},$ we have that all
the sets are not neighboring, and, as above, the sum of the corresponding
double probabilities is negligibly small with respect to any single sum over
the partitions of $\mathcal{K}_{0i},$ $i=1,...,n.$ Furthermore, the
probability $P(\mathcal{K}_{0}\cap\mathbb{K}_{\alpha/2},u)$ is also infinitely
smaller than the single sums. Finally, as it was mentioned above, the only
summands with slowest order give contribution in the final asymptotic
behavior, we may take the multipliers $1+o(1)$ out of the sum. Thus
Proposition is established.

\subsection{Main result}

Now we collect all obtained above asymptotic relations.

\begin{theorem}
\label{main}Let $S$ be a bounded open set in $\mathbb{R}^{d}$ containing zero,
and $X(\mathbf{t),}$ $\mathbf{t\in}S,$ be an a.s. continuous zeromean Gaussian
field satisfying Conditions \ref{CondA1} - \ref{CondA6}. Then for the
probability $P(S,u)$ given by (\ref{P(S,u)}) the following asymptotic
relations take places as $u\rightarrow\infty.$

\begin{itemize}
\item If $\dim\mathcal{K}_{0}=d,$ $P(S;u)$ satisfies the relation
(\ref{case S}).

\item If $\dim\mathcal{K}_{0}\in\lbrack1,d-1],$ $P(S;u)=P(\mathcal{K}%
_{0};u)(1+o(1)),$ and $\mathcal{K}_{0}$ satisfies the relation (\ref{K0_fin2}).

\item If $\dim\mathcal{K}_{0}=0$ and $\dim K_{c}>0,$ $P(S;u)=P(\mathcal{K}%
_{c};u)(1+o(1)),$ and $P(\mathcal{K}_{c};u)$ satisfies the relation
(\ref{K_c_fin}).

\item If $\dim\mathcal{K}_{0}=\dim K_{c}=0,$ $P(S;u)$ satisfies the relation
(\ref{case T}).
\end{itemize}
\end{theorem}

\textbf{Proof: }Since $B_{u}\subset\mathcal{K}_{0}\cup K_{c}\cup K_{\infty}$
and both $P(\mathcal{K}_{c};u)$ and $P(\mathcal{K}_{\infty};u)$ have order
$\Psi(u)$ which is infinitely smaller than $P(\mathcal{K}_{0};u)$ provided
$\mathcal{K}_{0}$ is not empty, the first two assertions follow from
Propositions \ref{stationary like} and \ref{Th_dimK_0<d}. The third assertion
follows from Propositions \ref{K_infty} and \ref{PP_gen}. The last relation is
given in Proposition \ref{Talagrand_case}.

Remark that the assertions of this Theorem agree with assertions of Theorem 3,
\cite{KHP}, where $d=1.$ There the case $\dim\mathcal{K}_{0}=1$ is described
in the items one ($\mathcal{K}_{0}$ consists of two manifolds) and two
($\mathcal{K}_{0}$ consists of one manifold, with the list of corresponding
cases for dimensions of $K_{c}$ and $K_{\infty}$). The case $\dim
\mathcal{K}_{0}=0$ is considered in the remainding items, with various
relations between dimensions of $K_{c}$ and $K_{\infty}.$

\section{Examples and discussion.}

First give examples of covariance functions satisfying the above conditions.

\subsection{\textbf{ }Example 1. Covariance functions of Pickands type.}

As we have seen, in one dimension case, $d=1,$ the only behavior of $1-r(t)$
satisfying Condition \ref{CondA4} is as following,
\begin{equation}
1-r(t)=|t|^{\alpha}\ell(|t|)(1+o(1)),\ \ t\rightarrow0, \label{Pic_d_1}%
\end{equation}
up to time scaling, with $\alpha\in(0,2]$ and slowly varying $\ell.$ In two
dimension case one can see two types of the behavior,
\begin{align}
1-r(\mathbf{t})  &  =(1+o(1))\sum_{i=1}^{2}|t_{i}|^{\alpha_{i}}\ell_{i}%
(|t_{i}|),\ \ \ \text{and}\nonumber\\
1-r(\mathbf{t})  &  =|\mathbf{t}|^{\alpha}\ell(|\mathbf{t}%
|)(1+o(1)),\ \ \mathbf{t\rightarrow0}, \label{Pic_d_2}%
\end{align}
up to linear time transformation of $\mathbb{R}^{2}$ and with the same
properties of $\alpha$s and $\ell$s. Observe that the first one can be a
covariance function of two independent stationary processes $\xi_{1}%
(t_{1})+\xi_{2}(t_{2}),$ with covariance functions as (\ref{Pic_d_1}). Going
this way, that is, summing independent fields with different time parameters,
one comes in $d$-dimension case to a \textquotedblleft
structured\textquotedblright covariance function,
\begin{equation}
1-r(\mathbf{t)=}(1+o(1))\sum_{j=1}^{n}|\mathbf{t}^{j}\mathbf{|}^{\alpha_{j}%
}\ell_{j}(|\mathbf{t}^{j}|),\ \mathbf{t\rightarrow0}, \label{structuredCF}%
\end{equation}
where $\mathbf{t}^{j}=(t_{l_{j}+1},...,t_{l_{j+1}}),$ $j=1,...,n,$ is a
partition of coordinates of $\mathbf{t}$, $0=l_{1}<...<l_{n}=d,$ $\alpha
_{j}\in(0,2],$ $j=1,...,n$. In the case without $\ell_{j}$s such functions are
considered in \cite{book}, the corresponding normalization $\mathbf{q}(u)$ is
also given there, namely
\begin{equation}
q_{i}(u)=u^{-2/\alpha_{j}},\ i=l_{j}+1,...,l_{j+1},\ j=0,...,n-1,
\label{structured_q(u)}%
\end{equation}
with
\begin{equation}
h(\mathbf{t})=\sum_{i=1}^{n}|\mathbf{t}_{i}\mathbf{|}^{\alpha_{i}}.
\label{structured_h(t)}%
\end{equation}
One can continue considering sums of independent fields, even with different
linear transforms of time parameters. A general case with $d=2$ is considered
in \cite{HDL}, see the following example.

\subsection{Example 2. The field from \cite{HDL}.}

In \cite{HDL} the covariance function satisfying%
\[
1-r(\mathbf{t})=\rho_{1}^{2}(|a_{11}t_{1}+a_{12}t_{2}|)+\rho_{2}^{2}%
(|a_{21}t_{1}+a_{22}t_{2}|),\ \mathbf{t\rightarrow0,}%
\]
is considered, where $\rho_{i},i=1,2,$ are regularly varying at zero functions
with indexes $\alpha_{i}\in(0,1].$ The rank of matrix $A=\left(
a_{ij},i,j=1,2\right)  $ can be $2,$ $1,$ or $0.$ The first case is described
here by Condition \ref{CondA4}. When the rank is equal to $1,$ the
standardized field in corresponding basis is equal to $X_{0}(t_{1}%
)\sigma(t_{1},t_{2}),$ with a Gaussian stationary process $X_{0}(t_{1}),$ that
is the field is degenerated along $t_{2}.$ In case of zero rank, the
standardized field is $X\sigma(t_{1},t_{2})$ with a Gaussian random variable
$X.$ Here such degenerated cases are not considered, but in a corresponding
basis one can represent $\mathbb{R}^{d}$ as a product of two spaces, dimension
of one of them should be equal to the rank of a matrix which generalized $A$
to $d$-dimension case, with subsequent application of given here results.

Notice that in \cite{HDL} the function $1-\sigma(\mathbf{t})$ is also assumed
to have a similar to $1-r(\mathbf{t})$ form, with some other regularly varying
functions. The corresponding matrix must be, of course, not degenerated,
otherwise one has infinitely many maximum points. From results here it follows
that such restriction on the variance is not necessary, in contrast of the
above representation for the covariance function. We would like to mention
that discussions with authors of \cite{HDL} helped us a lot in formulation of
our Conditions.

\subsection{On behavior of the variance.}

An example of very gentle behavior of the variance function at zero is
considered in \cite{KHP}, Example 1. Its trivial $d$-dimension generalization
can be as following,%
\[
1-\sigma^{2}(\mathbf{t})=\exp\left(  -|\mathbf{t}|^{-\beta(\mathbf{t/|t|)}%
}\right)  (1+o(1))
\]
as $\mathbf{t\rightarrow0,}$ where $\beta(\mathbf{e)}$ is a strictly positive
function given on the unit sphere $\mathbb{S}_{d-1}.$ Such behavior is a
subject of Proposition \ref{stationary like}, the behavior of the integral
$L_{f}(\lambda)$ as $\lambda\rightarrow\infty$ can be investigated similarly
to that in \cite{KHP}.

A generalization of the model from \cite{HDL} is
\[
1-\sigma^{2}(\mathbf{t})=|\mathbf{t}|^{\beta(\mathbf{t/|t|)}}\ell
_{\mathbf{t/|t|}}(|\mathbf{t}|)(1+o(1)),
\]
with similar $\beta(\mathbf{e})$ and a collection of slowly variable at zero
functions $\ell_{\mathbf{e}}(t),$ $\mathbf{e\in}\mathbb{S}_{d-1}.$ Remark
again that, in contrast of Pickands' behavior of the covariance function at
zero, the behavior of $\beta(\mathbf{e)}$ can be very variable.

An example when in the stationary like case $1-\sigma^{2}(\mathbf{t})$ is
close to $1-r(\mathbf{s,t}),$ $\mathbf{s,t\rightarrow0},$ is also considered
in \cite{KHP} for the one dimension case. Let $\ell_{\mathbf{e}}(t)$,
$|\mathbf{e}|=1$, be a family of slowly varying at zero functions. Take
$1-\sigma^{2}(\mathbf{t})=(1-r(\mathbf{t}))\ell_{\mathbf{t/|t|}}%
(|\mathbf{t}|)$ and denote $\mathbf{t}_{u}:=\mathbf{q}(u)\mathbf{t.}$When
$\ell_{\mathbf{t}_{u}\mathbf{/|t}_{u}\mathbf{|}}(|\mathbf{t}_{u}%
|)\rightarrow0$ as $u\rightarrow\infty,$ we have $\mathbf{t\in}\mathcal{K}%
_{0},$ when the limit is equal to infinity, $\mathbf{t\in}\mathcal{K}_{\infty
},$ so on.

One should assume, following Condition \ref{CondA6}, that the limit changes
from zero to non-zero and to infinity at most finite number of times. In
\cite{KHP}, $\ell(t)=\log(1/|t|)$ is taken, that is, only the stationary like
case is considered.

It is possible to evaluate similarly to \cite{KHP} asymptotic behavior of
$P(S;u)$ for particular cases of $R(\mathbf{s,t})$ and described here behavior
of $\sigma(\mathbf{t})$

\end{document}